\documentclass[12pt,psfig]{article}
\usepackage{graphicx,epsfig}
\usepackage{amsmath}
\usepackage{mathrsfs}
\usepackage{amssymb}
\def\ge{\geq}

\def\bbb{\begin{eqnarray*}}

\def\eee{\end{eqnarray*}}

\begin{document}

\baselineskip=17pt
\begin{center}

\vspace{-0.6in} {\large \bf Some criteria of chaos in non-autonomous discrete systems}\\ [0.2in]

Hua Shao$^\dag$, Guanrong Chen$^\dag$, Yuming Shi$^\ddag$

\vspace{0.15in}$^\dag$ Department of Electronic Engineering, City University of Hong Kong,\\
Hong Kong SAR, P.~R. China\\

\vspace{0.15in}$^\ddag$ Department of Mathematics, Shandong University \\
 Jinan, Shandong 250100, P.~R. China\\
\end{center}

{\bf Abstract.} This paper establishes some criteria of chaos in non-autonomous discrete systems.
Several criteria of strong Li-Yorke chaos are given. Based on these results, some criteria of
distributional chaos in a sequence are established. Moreover, several criteria of distributional
chaos induced by coupled-expansion for an irreducible transition matrix are obtained. Some of these
results not only extend the existing related results for autonomous discrete systems to non-autonomous
discrete systems, but also relax the assumptions of the counterparts. One example is provided for illustration.\medskip

\medskip

{2010 {\bf \it Mathematics Subject Classification}}: 37B55, 37D45, 37B10.

\bigskip

\noindent{\bf 1. Introduction}\medskip

Chaos of the non-autonomous discrete system (briefly, NDS)
\vspace{-0.2cm}$$x_{n+1} = f_n(x_n), \;\;n\geq0,                                                                                     \eqno(1.1)\vspace{-0.2cm}$$
has attracted considerable attention recently
\cite{Canovas,Huangqiu,kolyada,Shao15,shaocsf,shao18,Shi09,shi12,Tian,wu,Zhangli,Zhuyu,Zhu16},
where $\{f_n\}_{n=0}^{\infty}$ is a sequence of maps from $X$ to $X$, with $(X,d)$ being a metric
space. Many complex systems of real-world problems in the fields of biology, physics, chemistry and
engineering, are indeed non-autonomous, putting the model (1.1) into focus. The positive orbit
$\{x_n\}^{\infty}_{n=0}$ of system {\rm(1.1)} starting from an initial point $x_0\in X$
is given by $x_n=f^{n}_{0}(x_0)$, where $f^n_0:=f_{n-1}\circ\cdots\circ f_0$, $n\geq1$.

Coupled-expansion is always associated with chaos
\cite{Canovas,ju16,ju17,Kennedy,Kim,Kulczycki,Shao15,Shi04,Shi09,shicsf,shi12,Yang,Zhangli,zhangijbc,zhang}.
In 2009, the concept of coupled-expansion for a transition matrix was extended to NDSs,
for which a criterion of strong Li-Yorke chaos was established in \cite{Shi09}. Later,
the assumptions in this result were weakened in \cite{Shao15}, where some new criteria
of strong Li-Yorke chaos were established via coupled-expansion for an irreducible transition
matrix in bounded and closed subsets of a complete metric space, and in compact subsets of
a metric space, respectively. Then, a new concept of weak coupled-expansion for a transition
matrix was introduced for NDSs, and several criteria of chaos induced by weak coupled-expansion
for an irreducible transition matrix were established in \cite{Zhangli}. Recently, Li-Yorke
$\delta$-chaos for some $\delta>0$ and distributional $\delta'$-chaos in a sequence for some $\delta'>0$
were proved to be equivalent for system {\rm(1.1)} in the case that $(X,d)$ is a compact metric
space and $f_n$ are continuous maps for all $n\geq0$ \cite{shaocsf}. Consequently, those criteria of strong
Li-Yorke chaos can be regarded as criteria of distributional chaos in a sequence for system
{\rm(1.1)}. It arouses our interest to investigate the relationship of coupled-expansion for a transition
matrix and distributional chaos in a sequence for system {\rm(1.1)}.
In the present paper, we establish some new criteria of strong Li-Yorke chaos for system {\rm(1.1)},
which not only relax the assumptions of the counterparts in the literature, but also can be regarded as
criteria of distributional chaos in a sequence under certain conditions.

Recall that distributional chaos is a special kind of distributional chaos in a sequence in the case that the
sequence is the set of all nonnegative integers. However, generally it is stronger than distributional chaos
in a sequence. In 2015, one criterion of distributional chaos induced by coupled-expansion
for a transition matrix was established for autonomous discrete systems \cite{Kim}. Motivated by this result,
here we are interested in finding some criteria of distributional chaos induced by weak coupled-expansion for a transition
matrix for the non-autonomous system {\rm(1.1)}.

The rest of the present paper is organized as follows. Section 2 presents some basic concepts and useful
lemmas. In Section 3, two criteria of strong Li-Yorke chaos for system {\rm(1.1)} are obtained,
which relax the assumptions of the counterparts in the literature (see Remark 3.1). Moreover, other
two criteria of strong Li-Yorke chaos are derived by applying the relationship of Li-Yorke chaos between
system {\rm(1.1)} and its induced system. In Section 4, some criteria of distributional chaos in a sequence
are established for system {\rm(1.1)} by the results obtained in Section 3. Furthermore, several criteria of
distributional chaos induced by weak coupled-expansion for an irreducible transition matrix are obtained.
Finally, an example is provided in Section 5 for illustration.

\bigskip

\noindent{\bf 2. Preliminaries}\medskip

In this section, some basic concepts and useful lemmas are presented.\medskip

For convenience, denote $f^n_i:=f_{i+n-1}\circ\cdots\circ f_i$ and $f^{-n}_i:=(f_{i}^n)^{-1}$, $i\geq0$,
$n\geq1$. Let $A,B$ be nonempty subsets of $X$. The boundary of $A$ is denoted by $\partial A$; the diameter
of $A$ is denoted by $d(A):=\sup\{d(x,y): x,y\in A\}$; and the distance between $A$ and $B$ is denoted by $d(A,B):=
\inf\{d(a,b): a\in A,b\in B\}$. The set of all nonnegative integers and positive integers are denoted by $\mathbf{N}$
and $\mathbf{Z^{+}}$, respectively.\medskip

\noindent{\bf Definition 2.1} (\cite{Shi09}, Definition 2.7).
System {\rm(1.1)} is said to be Li-Yorke $\delta$-chaotic for some $\delta>0$ if it has an uncountable Li-Yorke
$\delta$-scrambled set $S$ in $X$; that is, for any $x, y\in S\subset X$,
\vspace{-0.2cm}$$\liminf_{n\to\infty}d(f_{0}^{n}(x),f_{0}^{n}(y))=0\;{\rm and}\;\limsup_{n\to\infty}d(f_{0}^{n}(x),f_{0}^{n}(y))>\delta.\vspace{-0.2cm}$$
Further, it is said to be chaotic in the strong sense of Li-Yorke if all the orbits starting from the points in $S$ are bounded.\medskip

\noindent{\bf Definition 2.2} (\cite{shaocsf}, Definitions 2.1 and 2.2).
System {\rm(1.1)} is said to be distributionally $\delta$-chaotic in a sequence $P=\{p_n\}_{n=1}^{\infty}$ for some $\delta>0$
if it has an uncountable distributionally $\delta$-scrambled set $D\subset X$ in $P$; that is, for any $x, y\in D\subset X$ and any $\epsilon>0$,
\vspace{-0.2cm}$$\limsup_{n\to\infty}\frac{1}{n}\sum_{i=1}^{n}\chi_{[0,\epsilon)}\big(d(f_{0}^{p_i}(x),f_{0}^{p_i}(y))\big)=1\;{\rm and}\;
\liminf_{n\to\infty}\frac{1}{n}\sum_{i=1}^{n}\chi_{[0,\delta)}\big(d(f_{0}^{p_i}(x),f_{0}^{p_i}(y))\big)=0,\vspace{-0.2cm}$$
where $\chi_{[0,\epsilon)}$ is the characteristic function defined on the set $[0,\epsilon)$. Further, if $P=\mathbf{N}$,
then system {\rm(1.1)} is said to be distributionally $\delta$-chaotic.\medskip

The relationship between Li-Yorke chaos and distributional chaos in a sequence for system {\rm(1.1)} is shown below.\medskip

\noindent{\bf Lemma 2.1} (\cite{shaocsf}, Theorem 3.6). {\it Let $(X, d)$ be compact and $f_n$ be continuous in $X$, $n\geq0$.
Then, system {\rm(1.1)} is Li-Yorke $\delta$-chaotic for some $\delta>0$ if and only if it is distributionally $\delta'$-chaotic
in a sequence for some $\delta'>0$.}\medskip

A matrix $A=(a_{ij})_{N\times N}$ ($N\geq2$) is said to be a transition matrix if $a_{ij}=0$ or $1$
for all $i,j$; $\sum_{j=1}^{N}a_{ij}\geq1$ for all $i$; and $\sum_{i=1}^{N}a_{ij}\geq1$ for all $j$,
$1\leq i,j\leq N$. A transition matrix $A=(a_{ij})_{N\times N}$ is said to be irreducible if, for
each pair $1\leq i, j\leq N$, there exists $k\in \mathbf{Z^{+}}$ such that $a^{(k)}_{ij}>0$, where
$a^{(k)}_{ij}$ denotes the $(i, j)$ entry of matrix $A^{k}$. A finite sequence $\omega=(s_1, s_2, \cdots, s_k)$
is said to be an allowable word of length $k$ for $A$ if $a_{s_{i}s_{i+1}}=1$, $1\leq i\leq k-1$, where $1\leq s_i\leq N$,
$1\leq i\leq k$. For convenience, the length of $\omega$ is denoted by $|\omega|$ and
$\omega^{n}:=\underbrace{(\omega,\cdots,\omega)}_{n}$, $n\geq1$.

The one-sided sequence space $\Sigma_{N}^{+}:=\{\alpha=(a_0,a_1,\cdots): 1\leq a_{i}\leq N,\;i\geq0\}$
is a metric space with the metric $\rho(\alpha, \beta):=\sum_{i=0}^{\infty}d'(a_i,b_i)/2^{i}$,
where $\alpha=(a_0,a_1,\cdots),\beta=(b_0,b_1,\cdots)\in\Sigma_{N}^{+}$, $d'(a_i, b_i)=0$ if $a_i=b_i$,
and $d'(a_i, b_i)=1$ if $a_i\neq b_i$, $i\geq0$. Note that $(\Sigma_{N}^{+},\rho)$ is a compact metric space.
Define the shift map $\sigma: \Sigma_{N}^{+}\to\Sigma_{N}^{+}$ by $\sigma((s_0,s_1,s_2,\cdots)):=(s_1,s_2,\cdots)$.
This map is continuous and $(\Sigma_{N}^{+},\sigma)$ is called the one-sided symbolic dynamical
system on $N$ symbols. For a given transition matrix $A=(a_{ij})_{N\times N}$, denote
\vspace{-0.2cm}$$\Sigma_{N}^{+}(A):=\{s=(s_0,s_1,\cdots): 1\leq s_{j}\leq N, \;a_{s_{j}s_{j+1}}=1, \;j\geq0\}.\vspace{-0.2cm}$$
$\Sigma_{N}^{+}(A)$ is a compact subset of $\Sigma_{N}^{+}$ and invariant under $\sigma$. The map
$\sigma_A:=\sigma|_{\Sigma_{N}^{+}(A)}: \Sigma_{N}^{+}(A)\to\Sigma_{N}^{+}(A)$ is said to be a
subshift of finite type for matrix $A$. For more details about symbolic dynamical systems
and subshifts of finite type, see \cite{Robinson,Zhou}.\medskip

The following two lemmas will also be useful in the sequel.\medskip

\noindent{\bf Lemma 2.2} (\cite{wangli}, Lemma 2.2). {\it $\Sigma_{2}^{+}$ has an uncountable subset $E$ such that
for any different points $\alpha=(a_0,a_1,\cdots),\beta=(b_0,b_1,\cdots)\in E$, $a_n=b_n$ for infinitely many
$n$ and $a_m\neq b_m$ for infinitely many $m$.}\medskip

\noindent{\bf Lemma 2.3} (\cite{zhangijbc}, Theorem 2.2). {\it Let $A=(a_{ij})_{N\times N}$ be an irreducible transition matrix with
$\sum_{j=1}^{N}a_{i_{0}j}\geq2$ for some $1\leq i_0\leq N$. Then
\begin{itemize}\vspace{-0.2cm}
\item[{\rm (i)}] for any $1\leq i, j\leq N$ and any $M\in\mathbf{Z^{+}}$, there exists at least one allowable word
$\omega=(i,\cdots,j)$ for $A$ such that $|\omega|>M$;\vspace{-0.2cm}
\item[{\rm (ii)}] for any given allowable word $\omega=(b_1, b_2, \cdots, b_k)$ for $A$, if $k>2N(N^{2}-2N+2)$, then
there exists another different allowable word $\omega'=(c_1, c_2,\cdots, c_k)$ for $A$ with $c_1=b_1$ and $c_k=b_k$.\vspace{-0.2cm}
\end{itemize}}\medskip

Next, the definition of weak coupled-expansion for a transition matrix is introduced.\medskip

\noindent{\bf Definition 2.3.} Let $A=(a_{ij})_{N\times N}$ be a transition matrix. If there exists a sequence
$\{V_{i,n}\}_{n=0}^{\infty}$ of nonempty subsets of $X$ with $V_{i,n}\cap V_{j,n}=\partial V_{i,n}\cap\partial V_{j,n}$
{\rm(}$d(V_{i,n},V_{j,n})>0${\rm)} for any $1\leq i\neq j\leq N$ and $n\geq0$ such that
\vspace{-0.2cm}$$f_n(V_{i,n})\supset\bigcup_{a_{ij}=1}V_{j,n+1},\;1\leq i\leq N, \;n\geq0,\vspace{-0.2cm}$$
then system {\rm(1.1)} is said to be {\rm(}strictly{\rm)} weakly $A$-coupled-expanding in $\{V_{i,n}\}_{n=0}^{\infty}$,
$1\leq i\leq N$. In the special case that $V_{i,n}=V_i$ for all $n\geq0$ and $1\leq i\leq N$, it is
said to be {\rm(}strictly{\rm)} $A$-coupled-expanding in $V_i$, $1\leq i\leq N$.\medskip

\noindent{\bf Remark 2.1.} Definition 2.3 is a slight revision of Definition 2.4 in \cite{Zhangli}.

\bigskip

\noindent{\bf 3. Some criteria of strong Li-Yorke chaos}\medskip

In this section, two criteria of strong Li-Yorke chaos for system {\rm(1.1)} are established.
Further, other two criteria of strong Li-Yorke chaos are given by applying the relationship
of Li-Yorke chaos between system {\rm(1.1)} and its induced system.\medskip

The following result can be easily verified.\medskip

\noindent{\bf Lemma 3.1.} {\it Let $A=(a_{ij})_{N\times N}$ be a transition matrix and $V_i$, $1\leq i\leq N$,
be disjoint nonempty compact subsets of $X$. Assume that $f_n$ is continuous in $\bigcup_{i=1}^{N}V_i$,
$n\geq0$, and $\{V_{i,n}\}_{n=0}^{\infty}$ is a sequence of nonempty closed subsets of $V_i$, $1\leq i\leq N$.
Then, $V_{\alpha}^{m,n}$ is compact and satisfies that $V_{\alpha}^{m+1,n}\subset V_{\alpha}^{m,n}\subset\bigcup_{i=1}^{N}
V_{i,n}$ for all $m,n\geq0$ and all $\alpha=(a_0,a_1,\cdots)$ $\in\Sigma_{N}^{+}(A)$, where
\vspace{-0.2cm}$$V_{\alpha}^{m,n}:=\bigcap_{k=0}^{m}f_{n}^{-k}(V_{a_{k},n+k}).                                                 \eqno(3.1)\vspace{-0.2cm}$$
Further, $(\bigcap_{m=0}^{\infty}V_{\alpha}^{m,n})\cap(\bigcap_{m=0}^{\infty}V_{\beta}^{m,n})=\emptyset$
for any $\alpha\neq\beta\in\Sigma_{N}^{+}(A)$ and $m,n\geq0$.}\medskip

\noindent{\bf Theorem 3.1.} {\it Let all the assumptions of Lemma 3.1 hold and suppose that $A$ is irreducible with
$\sum_{j=1}^{N}a_{i_{0}j}\geq2$ for some $1\leq i_0\leq N$. If
\begin{itemize}\vspace{-0.2cm}
\item[{\rm (i)}] $V_{\alpha}^{m,n}\neq\emptyset$ for all $m,n\geq0$ and all $\alpha=(a_0,a_1,\cdots)\in\Sigma_{N}^{+}(A)$,
where $V_{\alpha}^{m,n}$ is specified in {\rm(3.1)};\vspace{-0.2cm}
\item[{\rm (ii)}] there exist $\beta=(b_0,b_1,\cdots)\in\Sigma_{N}^{+}(A)$, an increasing subsequence
$\{m_k\}_{k=1}^{\infty}\subset\mathbf{Z^{+}}$, and $1\leq s_0\leq N$, such that $b_{m_{k}}=s_0$,
$k\geq1$, and $d(V_{\beta}^{m_k,n_k})$ converges to $0$ as $k\to\infty$ for all increasing subsequence
$\{n_k\}_{k=1}^{\infty}\subset\mathbf{N}$;\vspace{-0.2cm}
\end{itemize}
then system {\rm(1.1)} is chaotic in the strong sense of Li-Yorke.}\medskip

\noindent{\bf Proof.} Since $A$ is a transition matrix, there exist $1\leq t_0,r_0\leq N$ such that $a_{t_0b_0}=a_{b_0r_0}=1$.
By Lemma 2.3, there exist four allowable words for matrix $A$:
\vspace{-0.2cm}$$\omega_0=(r_0,\cdots,t_0),\;\omega_1=(b_0,\cdots,t_0),\;\omega_2=(b_0,\cdots,t_0),\;\omega_3=(s_0,\cdots,b_0),\vspace{-0.2cm}$$
where $|\omega_0|=l_1$, $|\omega_1|=|\omega_2|=l_2$ with $\omega_1\neq\omega_2$, and $|\omega_3|=l_3$
(if $s_0=b_0$, then set $\omega_3=(b_0)$ with length $1$). Denote
\vspace{-0.2cm}$$\Omega:=\{\gamma=(B_1, B_2,\cdots): B_i=\omega_1\; {\rm or}\; \omega_2,\; i\geq1\}.                           \eqno(3.2)\vspace{-0.2cm}$$
Then $\Omega\subset\Sigma_{N}^{+}(A)$ is uncountable. Note that $a_{b_{m_{k}-1}s_0}=1$, $k\geq1$, since $b_{m_{k}}=s_0$.
For any $\gamma=(B_1, B_2,\cdots)\in\Omega$, set
\vspace{-0.2cm}$$\hat{\gamma}:=(b_{0},\cdots,b_{m_{1}-1},\omega_3,\omega_0,B_1,b_{0},\cdots,b_{m_{2}-1},\omega_3,\omega_0,B_1,B_2,\cdots). \eqno(3.3)\vspace{-0.2cm}$$
Clearly, $\hat{\gamma}\in\Sigma_{N}^{+}(A)$ and $\hat{\gamma_1}\neq\hat{\gamma_2}$ if and only if $\gamma_1\neq\gamma_1$.
By assumption (i) and Lemma 3.1, one has that $V_{\hat{\gamma}}^{m,0}$ is nonempty, compact, and satisfies that
$V_{\hat{\gamma}}^{m+1,0}\subset V_{\hat{\gamma}}^{m,0}$, $m\geq0$. Thus
$\bigcap_{m=0}^{\infty}V_{\hat{\gamma}}^{m,0}\neq\emptyset$.
Fix one point $x_{\hat{\gamma}}\in\bigcap_{m=0}^{\infty}V_{\hat{\gamma}}^{m,0}$, and denote
\vspace{-0.2cm}$$S:=\{x_{\hat{\gamma}}: \gamma\in\Omega\}.\vspace{-0.2cm}$$
It follows from Lemma 3.1 that $x_{\hat{\gamma_1}}\neq x_{\hat{\gamma_2}}$ if and only if
$\hat{\gamma_1}\neq\hat{\gamma_2}$, which holds if and only if $\gamma_1\neq\gamma_2$.
So, $S$ is uncountable.

Next, it will be shown that $S$ is a Li-Yorke $\delta$-scrambled set for system {\rm(1.1)}, where
$\delta:=\min_{1\leq i\neq j\leq N}d(V_i, V_j)>0$. For any $x_{\hat{\gamma_1}}\neq x_{\hat{\gamma_2}}\in S$,
one has that $\hat{\gamma_1}\neq\hat{\gamma_2}$. By (3.1)-(3.3), there exist an infinite sequence
$\{h_k\}_{k=1}^{\infty}$ and $1\leq s_1\neq s_2\leq N$ such that $f_{0}^{h_k}(x_{\hat{\gamma_1}})\in V_{s_1,h_k}\subset V_{s_1}$
and $f_{0}^{h_k}(x_{\hat{\gamma_2}})\in V_{s_2,h_k}\subset V_{s_2}$, $k\geq1$.
So,
\vspace{-0.2cm}$$d(f_{0}^{h_k}(x_{\hat{\gamma_1}}),f_{0}^{h_k}(x_{\hat{\gamma_2}}))\geq\delta,\;\;k\geq1.\vspace{-0.2cm}$$
Thus,
\vspace{-0.2cm}$$\limsup_{n\to\infty}d(f_{0}^{n}(x_{\hat{\gamma_1}}),f_{0}^{n}(x_{\hat{\gamma_2}}))\geq\delta.                 \eqno(3.4)\vspace{-0.2cm}$$
On the other hand, by (3.3), one has that $\sigma_{A}^{n_{k}}(\hat{\gamma})=(b_{0},\cdots,b_{m_{k}},\cdots)$ for any $\gamma\in\Omega$, where
\[n_k=(k-1)l_1+\frac{k(k-1)}{2}l_2+(k-1)l_3+\sum_{t=1}^{k}m_t,\;k\geq1.\]
Then, $f_{0}^{n_{k}+i}(x_{\hat{\gamma_1}}),f_{0}^{n_{k}+i}(x_{\hat{\gamma_2}})\in V_{b_i,n_k+i}$, $0\leq i\leq m_k$.
Thus, $f_{0}^{n_{k}}(x_{\hat{\gamma_1}}),f_{0}^{n_{k}}(x_{\hat{\gamma_2}})\in V_{\beta}^{m_{k},n_{k}}$, $k\geq1$,
implying that
\vspace{-0.2cm}$$d(f_{0}^{n_{k}}(x_{\hat{\gamma_1}}),f_{0}^{n_{k}}(x_{\hat{\gamma_2}}))\leq d(V_{\beta}^{m_{k},n_{k}}),\;k\geq1.\vspace{-0.2cm}$$
This, together with assumption (ii), yields that
\vspace{-0.2cm}$$\lim_{k\to\infty}d(f_{0}^{n_{k}}(x_{\hat{\gamma_1}}),f_{0}^{n_{k}}(x_{\hat{\gamma_2}}))=0.\vspace{-0.2cm}$$
So,
\vspace{-0.2cm}$$\liminf_{n\to\infty}d(f_{0}^{n}(x_{\hat{\gamma_1}}),f_{0}^{n}(x_{\hat{\gamma_2}}))=0.                               \eqno(3.5)\vspace{-0.2cm}$$
Hence, $S$ is an uncountable $\delta$-scrambled set of system {\rm(1.1)} by (3.4) and (3.5).
Moreover, for any $x_{\hat{\gamma}}\in S$, its positive orbit
$\{f_{0}^{n}(x_{\hat{\gamma}})\}_{n=0}^{\infty}\subset\bigcup_{i=1}^{N}V_{i}$ by (3.1) and (3.3).
Thus, all the orbits starting from the points in $S$ are bounded.
Therefore, system {\rm(1.1)} is chaotic in the strong sense of Li-Yorke.
The proof is complete.\medskip

\noindent{\bf Remark 3.1.} It can be easily verified that the weak $A$-coupled-expansion in $\{V_{i,n}\}_{n=0}^{\infty}$,
$1\leq i\leq N$, implies the condition in assumption {\rm(i)} of Theorem 3.1. However, the converse is not true in general,
even in the special case that $f_n=f$ and $V_{i,n}=V_i$, $n\geq0$, $1\leq i\leq N$ (see Example 3.1.1 in \cite{ju17}).
Hence, assumption {\rm(i)} of Theorem 3.1 here is strictly weaker than assumption ${\rm(ii_{1})}$ of Theorems 3.2 and 3.3
in \cite{Zhangli}. Moreover, both assumption ${\rm(iii_{1})}$ of Theorem 3.2 and assumption ${\rm(iii_{4})}$ of Theorem 3.3 in
\cite{Zhangli} imply assumption {\rm(ii)} of Theorem 3.1 above (see, for example, the proof of Corollary 4.1 below).
Hence, Theorems 3.2 and 3.3 in \cite{Zhangli} are corollaries of Theorem 3.1 here.\medskip

The following result is a direct consequence of Theorem 3.1.\medskip

\noindent{\bf Corollary 3.1.} {\it Let all the assumptions of Theorem 3.1 hold except that assumption {\rm(i)} is replaced by
that system {\rm(1.1)} is weakly $A$-coupled-expanding in $\{V_{i,n}\}_{n=0}^{\infty}$, $1\leq i\leq N$.
Then, system {\rm(1.1)} is chaotic in the strong sense of Li-Yorke.}\medskip

Let $\{k_n\}_{n=1}^{\infty}\subset\mathbf{Z^{+}}$ be an increasing subsequence. The following system:
\vspace{-0.2cm}$$\hat{x}_{n+1}=\hat{f}_n(\hat{x}_n),\;n\geq0,                                                                  \eqno(3.6)\vspace{-0.2cm}$$
is called the induced system by system {\rm(1.1)} through $\{k_n\}_{n=1}^{\infty}$ \cite{shi12}, where
\vspace{-0.2cm}$$\hat{f}_0:=f_{k_{1}-1}\circ f_{k_{1}-2}\circ\cdots f_0,\;\;\hat{f}_n:=f_{k_{n+1}-1}\circ f_{k_{n+1}-2}\circ\cdots f_{k_n},\;\;n\geq1.\vspace{-0.2cm}$$
Let $\{x_n\}_{n=0}^{\infty}$ be the orbit of system (1.1) starting from $x_0$ and $\{\hat{x}_n\}_{n=0}^{\infty}$
be the orbit of the induced system (3.6) starting from $\hat{x}_0:=x_0$.
Then, $\hat{x}_n=f_{0}^{k_n}(x_0)$, $n\geq0$, and thus the orbit $\{\hat{x}_n\}_{n=0}^{\infty}$
of system (3.6) is a part of the orbit of system {\rm(1.1)} starting from the same initial point $x_0$.\medskip

Next, recall the relationship of (strong) Li-Yorke chaos between systems {\rm(1.1)} and {\rm(3.6)}.\medskip

\noindent{\bf Lemma 3.2} (\cite{shi12}, Theorem 3.1). {\it If system {\rm(3.6)} is Li-Yorke $\delta$-chaotic for some $\delta>0$
through $\{k_n\}_{n=1}^{\infty}$, so is system {\rm(1.1)}. Further, if system {\rm(3.6)} is chaotic in the strong sense of Li-Yorke
through $\{k_n\}_{n=1}^{\infty}$, so is system {\rm(1.1)} in the case that $X$ is bounded.}\medskip

The following result follows directly from Theorem 3.1 and Lemma 3.2.\medskip

\noindent{\bf Theorem 3.2.} {\it Assume that there exists an increasing subsequence $\{k_n\}_{n=1}^{\infty}\subset\mathbf{Z^{+}}$
such that all the assumptions of Theorem 3.1 hold for system {\rm(3.6)}. Then system {\rm(1.1)} is Li-Yorke $\delta$-chaotic
for some $\delta>0$. Further, system {\rm(1.1)} is chaotic in the strong sense of Li-Yorke in the case that $X$ is bounded.}\medskip

Applying Lemma 3.2 to Corollary 3.1, one obtains the following result.\medskip

\noindent{\bf Corollary 3.2.} {\it Assume that there exists an increasing subsequence $\{k_n\}_{n=1}^{\infty}\subset\mathbf{Z^{+}}$
such that all the assumptions of Corollary 3.1 hold for system {\rm(3.6)}. Then, system {\rm(1.1)} is Li-Yorke $\delta$-chaotic for
some $\delta>0$. Further, system {\rm(1.1)} is chaotic in the strong sense of Li-Yorke in the case that $X$ is bounded.}\bigskip

\noindent{\bf 4. Some criteria of distributional chaos in a sequence and distributional chaos}\medskip

In this section, several criteria of distributional chaos in a sequence and of distributional chaos are established
for system {\rm(1.1)},  respectively.\medskip

Applying Lemma 2.1 to Theorems 3.1 and 3.2 and Corollaries 3.1 and 3.2, respectively, the following result can be obtained.\medskip

\noindent{\bf Theorem 4.1.} {\it Let all the assumptions of Theorem 3.1 or Corollary 3.1 or Theorem 3.2 or Corollary 3.2 hold.
If $X$ is compact and $f_n$ is continuous in $X$ for all $n\geq0$, then system {\rm(1.1)} is distributionally $\delta$-chaotic
in a sequence for some $\delta>0$.}\medskip

The next result gives a criterion of distributional chaos for system (1.1), which is induced by weak coupled-expansion
for an irreducible transition matrix.\medskip

\noindent{\bf Theorem 4.2.} {\it Let all the assumptions of Lemma 3.1 hold, $\{f_n\}_{n=0}^{\infty}$ be equi-continuous in
$\bigcup_{i=1}^{N}V_i$, and $A$ be irreducible with $\sum_{j=1}^{N}a_{i_{0}j}\geq2$ for some $1\leq i_0\leq N$. If
\begin{itemize}\vspace{-0.2cm}
\item[${\rm(i)}$] system {\rm(1.1)} is weakly $A$-coupled-expanding in $\{V_{i,n}\}_{n=0}^{\infty}$, $1\leq i\leq N$;\vspace{-0.2cm}
\item[${\rm(ii)}$] there exists a periodic point $\gamma=(r_0,r_1,\cdots)\in\Sigma_{N}^{+}(A)$ such that
$d(V_{\gamma}^{m,n})$ uniformly converges to $0$ with respect to $n\geq0$ as $m\to\infty$;\vspace{-0.2cm}
\end{itemize}
then system {\rm(1.1)} is distributionally $\delta$-chaotic for some $\delta>0$.}\medskip

\noindent{\bf Proof.} Since $A$ is a transition matrix, there exist $1\leq l_0,m_0\leq N$ such that $a_{l_0r_0}=a_{r_0m_0}=1$.
By Lemma 2.3, one has that there exist three allowable words for matrix $A$:
\vspace{-0.2cm}$$\omega_0:=(m_0,\cdots,l_0),\;\omega_1:=(r_0,\cdots,l_0),\;\omega_2:=(r_0,\cdots,l_0),\vspace{-0.2cm}$$
where $|\omega_1|=|\omega_2|=l$ and $\omega_1\neq\omega_2$.
Let $E\subset\Sigma_{2}^{+}$ be the set satisfying the property in Lemma 2.2.
For any $\beta=(b_0,b_1,b_2,\cdots)\in E$, set
\vspace{-0.2cm}$$\hat{\beta}:=(r_0,\omega_0,\omega_{b_0}^{p_1},r_0,\cdots,r_{m_{2}-2},R_2,\omega_0,
\omega_{b_1}^{p_2},r_0,\cdots,r_{m_{3}-2},R_3,\omega_0,\omega_{b_2}^{p_3},\cdots),\eqno(4.1)\vspace{-0.2cm}$$
where $p_1:=2|(r_0,\omega_0)|$, $p_n=2^{n}|(r_0,\omega_0,\cdots,R_n,$ $\omega_0)|$ for $n\geq2$;
$m_n:=2^{n}|(r_0,\omega_0,\cdots,\omega_{b_{n-2}}^{p_{n-1}})|$ for $n\geq2$;
and for any $n\geq2$, $R_n:=(r_{0})$ if $r_{m_{n}-1}=r_0$, and otherwise, $R_n:=(r_{m_{n}-1},\cdots,r_0)$,
while the fact that there exists an allowable word $(r_{m_{n}-1},\cdots,r_0)$ for the matrix $A$ since $A$
is irreducible has been used in the case that $r_{m_{n}-1}\ne r_0$. Clearly, $\hat{\beta}\in\Sigma_{N}^{+}(A)$
and $\hat{\beta_1}\neq\hat{\beta_2}$ if and only if $\beta_1\neq\beta_2$. By assumption (i) and Lemma 3.1,
$V_{\hat{\beta}}^{m,0}$ is nonempty, compact, and satisfies that $V_{\hat{\beta}}^{m+1,0}
\subset V_{\hat{\beta}}^{m,0}$, $m\geq0$. Thus, $\bigcap_{m=0}^{\infty}V_{\hat{\beta}}^{m,0}
\neq\emptyset$. Fix one point $x_{\hat{\beta}}\in\bigcap_{m=0}^{\infty}V_{\hat{\beta}}^{m,0}$, and denote
\vspace{-0.2cm}$$D:=\{x_{\hat{\beta}}: \beta\in E\}.\vspace{-0.2cm}$$
Using Lemma 3.1 again, one has that $x_{\hat{\beta_1}}\neq x_{\hat{\beta_2}}$ if and only if
$\hat{\beta_1}\neq\hat{\beta_2}$, which holds if and only if $\beta_1\neq\beta_2$.
Hence, $D$ is uncountable since $E$ is uncountable.

In the following, it will be shown that $D$ is a distributionally $\delta_0$-scrambled subset for system (1.1), where
$\delta_0:=\min_{1\leq i\neq j\leq N}d(V_i,V_j)$. Fix any $x_{\hat{\beta_1}},x_{\hat{\beta_2}}\in D$ with $x_{\hat{\beta_1}}\neq x_{\hat{\beta_2}}$.
Then, $\beta_1\neq\beta_2\in E$. Suppose that $\beta_1=(b_{0}^{(1)},b_1^{(1)},b_2^{(1)},\cdots)$ and $\beta_2=(b_0^{(2)},b_1^{(2)},b_2^{(2)},\cdots)$.
Thus, there exists an increasing subsequence $\{t_i\}_{i=1}^{\infty}\subset\mathbf{Z^{+}}$ such that $b_{t_i}^{(1)}
\neq b_{t_i}^{(2)}$, $i\geq1$, by Lemma 2.2. Set
\vspace{-0.2cm}$$n_i:=|(r_0,\omega_0,\cdots,\omega_{b_{t_i}}^{p_{{t_i}+1}})|,\;i\geq1.\vspace{-0.2cm}$$
Then,
\vspace{-0.2cm}$$n_i=p_{{t_i}+1}l+p_{{t_i}+1}2^{-(t_{i}+1)},\;i\geq1.                                                          \eqno(4.2)\vspace{-0.2cm}$$
Denote
\vspace{-0.2cm}$$\Omega_0:=\{\omega: \omega \;{\rm is\; an\; allowable\; word\; for}\; A\; {\rm with}\;\; |\omega|=l\}.\vspace{-0.2cm}$$
For any $\mathcal{C}=(c_1,\cdots,c_{l})\in\Omega_0$ and any $n\geq0$, set
\vspace{-0.2cm}$$V_{\mathcal{C}}^{n}:=\bigcap_{k=0}^{l-1}f_{n}^{-k}(V_{c_{k+1},n+k}).\vspace{-0.2cm}$$
It is evident that
\vspace{-0.2cm}$$d_n:=\inf\{d(V_{\mathcal{C}}^{n},V_{\mathcal{D}}^{n}): \mathcal{C}\neq\mathcal{D}\in\Omega\}\geq\delta_0,\;n\geq0.\eqno(4.3)\vspace{-0.2cm}$$
By (4.1) and (4.3), one has that
\vspace{-0.2cm}$$\sum_{j=0}^{n_{i}-1}\chi_{[0,\delta_0)}(d(f_{0}^{j}(x_{\hat{\beta_1}}),f_{0}^{j}(x_{\hat{\beta_2}})))\leq n_{i}-(p_{t_{i}+1}l-l+1),\;\;i\geq1.\vspace{-0.2cm}$$
This, together with (4.2), implies that
\vspace{-0.2cm}$$\lim_{i\to\infty}\frac{1}{n_i}\sum_{j=0}^{n_{i}-1}\chi_{[0,\delta_0)}\big(d(f_{0}^{j}(x_{\hat{\beta_1}}),f_{0}^{j}(x_{\hat{\beta_2}}))\big)
\leq1-\lim_{i\to\infty}\frac{p_{t_{i}+1}l-l+1}{p_{{t_i}+1}l+p_{{t_i}+1}2^{-(t_{i}+1)}}=0,\vspace{-0.2cm}$$
which yields that
\vspace{-0.2cm}$$\liminf_{n\to\infty}\frac{1}{n}\sum_{j=0}^{n-1}\chi_{[0,\delta_0)}\big(d(f_{0}^{j}(x_{\hat{\beta_1}}),f_{0}^{j}(x_{\hat{\beta_2}}))\big)=0.
\eqno(4.4)\vspace{-0.2cm}$$
On the other hand, denote
\vspace{-0.2cm}$$k_i:=|(r_0,\omega_0,\cdots,r_0,\cdots,r_{m_{i}-1})|,\;i\geq1.\vspace{-0.2cm}$$
Thus,
\vspace{-0.2cm}$$k_i=m_i+m_i2^{-i},\;i\geq1.                                                                                   \eqno(4.5)\vspace{-0.2cm}$$
Fix any $\epsilon>0$. It is claimed that there exists $M>0$ such that, for all $m\geq M$,
\vspace{-0.2cm}$$d(V_{\sigma_{A}^{j}(\gamma)}^{m,n+j})<\epsilon,\;0\leq j\leq P-1,\;n\geq0,                                    \eqno(4.6)\vspace{-0.2cm}$$
where $P$ is the period of $\gamma$. For simplicity, only consider the case of $P=2$.
The general cases can be proved in a similar way. Since $\{f_n\}_{n=0}^{\infty}$ is equi-continuous
in $\bigcup_{i=1}^{N}V_i$, there exists $0<\delta_1<\epsilon$ such that, for any $x,y\in X$ with $d(x,y)<\delta_1$,
\vspace{-0.2cm}$$d(f_n(x),f_n(y))<\epsilon,\;n\geq0.                                                                           \eqno(4.7)\vspace{-0.2cm}$$
By assumption (ii), there exists $M>0$ such that, for any $m\geq M$,
\vspace{-0.2cm}$$d(V_{\gamma}^{m,n})<\delta_1<\epsilon,\;n\geq0.                                                               \eqno(4.8)\vspace{-0.2cm}$$
Fix any $n\geq0$ and $m\geq M$. By assumption (i), one has that
\vspace{-0.2cm}$$f_n(V_{\gamma}^{m+1,n})=V_{\sigma_{A}(\gamma)}^{m,n+1}.                                                       \eqno(4.9)\vspace{-0.2cm}$$
It follows from (4.9) that, for any $z_1,z_2\in V_{\sigma_{A}(\gamma)}^{m,n+1}$, there exist $x,y\in V_{\gamma}^{m+1,n}$ such that
$z_1=f_n(x)$ and $z_2=f_n(y)$. By (4.8), one has that
\vspace{-0.2cm}$$d(x,y)\leq d(V_{\gamma}^{m+1,n})<\delta_1.\vspace{-0.2cm}$$
This, together with (4.7), implies that
\vspace{-0.2cm}$$d(z_1,z_2)=d(f_n(x),f_n(y))<\epsilon.\vspace{-0.2cm}$$
Thus,
\vspace{-0.2cm}$$d(V_{\sigma_{A}(\gamma)}^{m,n+1})<\epsilon.                                                                   \eqno(4.10)\vspace{-0.2cm}$$
Hence, (4.6) follows from (4.8) and (4.10) in the case of $P=2$.
By (4.1) and (4.6), one has that
\vspace{-0.2cm}$$\sum_{j=0}^{k_{i}-1}\chi_{[0,\epsilon)}\big(d(f_{0}^{j}(x_{\hat{\beta_1}}),f_{0}^{j}(x_{\hat{\beta_2}}))\big)\geq m_{i}-M,\;i\geq1.\vspace{-0.2cm}$$
This, together with (4.5), implies that
\vspace{-0.2cm}$$\lim_{i\to\infty}\frac{1}{k_i}\sum_{j=0}^{k_{i}-1}\chi_{[0,\epsilon)}\big(d(f_{0}^{j}(x_{\hat{\beta_1}}),f_{0}^{j}(x_{\hat{\beta_2}}))\big)
\geq\lim_{i\to\infty}\frac{m_{i}-M}{k_{i}}=\lim_{i\to\infty}\frac{m_{i}-M}{m_i+m_i2^{-i}}=1.\vspace{-0.2cm}$$
So,
\vspace{-0.2cm}$$\limsup_{n\to\infty}\frac{1}{n}\sum_{j=0}^{n-1}\chi_{[0,\epsilon)}\big(d(f_{0}^{j}(x_{\hat{\beta_1}}),f_{0}^{j}(x_{\hat{\beta_2}}))\big)=1,
\;\forall\;\epsilon>0.\eqno(4.11)\vspace{-0.2cm}$$

Hence, $D$ is an uncountable distributionally $\delta_0$-scrambled subset for system {\rm(1.1)} by (4.4) and (4.11).
Therefore, system {\rm(1.1)} is distributionally $\delta_0$-chaotic.
The proof is complete.\medskip

\noindent{\bf Remark 4.1.} Theorem 4.2 extends Theorem 3.2 in \cite{Kim} from autonomous to non-autonomous systems.\medskip

\noindent{\bf Corollary 4.1.} {\it Let all the assumptions of Theorem 4.2 hold, except that assumption {\rm(ii)}
is replaced by
\begin{itemize}
\item[${\rm(ii_a)}$] there exists a constant $\lambda>1$ such that
\vspace{-0.2cm}$$d(f_n(x),f_n(y))\geq \lambda d(x,y), \forall\;x, y\in V_{j,n},\;1\leq j\leq N,\;n\geq0,                          \eqno(4.11)\vspace{-0.2cm}$$
\end{itemize}
then system {\rm(1.1)} is distributionally $\delta$-chaotic for some $\delta>0$.}\medskip

\noindent{\bf Proof.} It can be shown that $d(V_{\gamma}^{m,n})$ uniformly converges to $0$ with respect to $n\geq0$
as $m\to\infty$ for all $\gamma=(r_0,r_1,\cdots)\in\Sigma_{N}^{+}(A)$. Indeed, it follows from (4.11) that,
for any $x,y\in V_{\gamma}^{m,n}$,
\vspace{-0.2cm}$$d(f_{n}^{m}(x), f_{n}^{m}(y))\geq \lambda^{m}d(x,y),\;m\geq1,\;n\geq0,\vspace{-0.2cm}$$
implying that
\vspace{-0.2cm}$$d(x,y)\leq\lambda^{-m}d(f_{n}^{m}(x),f_{n}^{m}(y))\leq\lambda^{-m}d(V_{j_0,n+m})\leq\lambda^{-m}d(V_{j_0}),\;m\geq1,\;n\geq0.\vspace{-0.2cm}$$
Thus,
\vspace{-0.2cm}$$d(V_{\gamma}^{m,n})\leq\lambda^{-m}d(V_{j_0}),\;m\geq1,\;n\geq0.\vspace{-0.2cm}$$
This, together with the assumption that $\lambda>1$, yields that $d(V_{\gamma}^{m,n})$ uniformly converges to $0$ with respect to
$n\geq0$ as $m\to\infty$. Hence, all the assumptions of Theorem 4.2 hold. Therefore, system {\rm(1.1)} is distributionally
$\delta$-chaotic for some $\delta>0$. This completes the proof.\medskip

The following result is somewhat better since it only requires $f_n$ be expanding in distance in one subset for all $n\geq0$.\medskip

\noindent{\bf Corollary 4.2.} {\it Let all the assumptions of Theorem 4.2 hold, except that assumption {\rm(ii)}
is replaced by
\begin{itemize}
\item[${\rm(ii_b)}$] there exist an integer $1\leq j_0\leq N$ and a constant $\lambda>1$ such that
$a_{j_{0}j_{0}}=1$ and
\vspace{-0.2cm}$$d(f_n(x),f_n(y))\geq \lambda d(x,y), \;\forall\;x, y\in V_{j_{0},n},\;n\geq0,\vspace{-0.2cm}$$
\end{itemize}
then system {\rm(1.1)} is distributionally $\delta$-chaotic for some $\delta>0$.}\medskip

\noindent{\bf Proof.} With a similar proof to that of Corollary 4.1, one can show that for $\gamma=(j_0,j_0,j_0,\cdots)$,
$d(V_{\gamma}^{m,n})$ uniformly converges to $0$ with respect to $n\geq0$ as $m\to\infty$.
Hence, all the assumptions of Theorem 4.2 hold. Therefore, system {\rm(1.1)} is distributionally
$\delta$-chaotic for some $\delta>0$. The proof is complete.

\bigskip

\noindent{\bf 5. An example}\medskip

In this section, an example is provided to illustrate the theoretical results.\medskip

\noindent{\bf Example 5.1.} Consider the following non-autonomous logistic system:
\vspace{-0.2cm}$$x_{n+1}=r_nx_n(1-x_n),\;n\geq0,                                                                               \eqno(5.1)\vspace{-0.2cm}$$
governed by the maps $f_n(x)=r_nx(1-x)$, $x\in[0,1]$, where $9/2\leq r_n\leq\mu$, $n\geq0$,
and $\mu\ge 9/2$ is a constant. It can be easily verified that
\vspace{-0.2cm}$$|f'_{n}(x)|=r_n|1-2x|\leq r_n\leq\mu,\;\forall\;x\in[0,1].\vspace{-0.2cm}$$
Thus,
\vspace{-0.2cm}$$|f_{n}(x)-f_{n}(y)|\leq\mu |x-y|,\;\forall\;x\in[0,1],\vspace{-0.2cm}$$
which yields that $\{f_n\}_{n=0}^{\infty}$ is equi-continuous in $[0,1]$. Set
\vspace{-0.2cm}$$V_1=[0,1/3],\;V_2=[3/5,1].\vspace{-0.2cm}$$
It is evident that $V_1$ and $V_2$ are disjoint nonempty compact subsets of $[0,1]$, and
\vspace{-0.2cm}$$f_n(V_1)=[0,6r_n/25],\;f_n(V_2)=[0,2r_n/9].\vspace{-0.2cm}$$
Since $r_n\geq9/2$, $n\geq0$, one can see that
\vspace{-0.2cm}$$V_1\cup V_2\subset f_n(V_1)\cap f_n(V_2),\;n\geq0,                                                            \eqno(5.2)\vspace{-0.2cm}$$
which implies that system (5.1) is strictly $A$-coupled-expanding in $V_1$ and $V_2$, where $A=(a_{ij})_{2\times2}$
with $a_{ij}=1$, $1\leq i,j\leq2$, and thus $A$ is irreducible with $\sum_{j=1}^{2}a_{i_{0}j}=2$, $i_0=1,2$.
On the other hand, one has that
\vspace{-0.2cm}$$|f'_{n}(x)|=r_n|1-2x|\geq \frac{r_n}{3}\geq\frac{3}{2},\;\forall\;x\in V_1,\vspace{-0.2cm}$$
which yields that
\vspace{-0.2cm}$$|f_{n}(x)-f_{n}(y)|\geq\frac{3}{2}|x-y|,\;\forall\;x,y\in V_1.                                                \eqno(5.3)\vspace{-0.2cm}$$

Therefore, all the assumptions in Theorem 3.1 hold for system {\rm(1.1)} with
$\beta=(1,1,\cdots)\in\Sigma_{2}^{+}(A)$, satisfying by (5.3) that $d(V_{\beta}^{m,n})$ uniformly converges to
$0$ with respect to $n\geq0$ as $m\to\infty$. By Theorem 3.1, system {\rm(1.1)} is chaotic
in the strong sense of Li-Yorke.

Moreover, all the assumptions in Corollary 4.2 hold for system {\rm(1.1)} with $j_0=1$, satisfying
assumption ${\rm(ii_b)}$ by (5.2) and (5.3). By Corollary 4.2, system {\rm(1.1)} is distributionally
$\delta$-chaotic for some $\delta>0$.\medskip

\noindent{\bf Remark 5.1.} System (5.1) is an important model in biology, which
describes the population growth under certain conditions. Comparing to Example 5.1 in \cite{Shi09},
here it not only proves that system (5.1) is chaotic in the strong sense of Li-Yorke, but also proves that
system (5.1) is distributionally $\delta$-chaotic for some $\delta>0$.

Finally, it is worth noting here that Theorem 4.2 and its corollaries still hold true when $f_n: X\to X$
is replaced by that $f_n: X_n\to X_{n+1}$ with compact subsets $X_n\subset X$ for any $n\geq0$.\medskip

\bigskip

\noindent{\bf Acknowledgments}\medskip

This research was supported by the Hong Kong Research Grants Council (GRF Grant CityU11200317) and the NNSF of China (Grant 11571202).

\end{document}